\documentclass[preprint,10pt]{elsarticle}

 \usepackage{esvect} 
 \usepackage{latexsym}
\usepackage{amsfonts}       % ,epsfig}
\usepackage{color}
\usepackage{amsmath}
\newcounter{eqnsave}

\def\C{{{{\rm {\mbox{\small l}}} \kern -.50em {\rm C}}}}
\def\I{{{{\rm l} \kern -.10em {\rm I}}}}
\def\R{{{{\rm l} \kern -.15em {\rm R}}}}
\def\N{{{{\rm l} \kern -.15em {\rm N}}}}
\def\E{{{{\rm l} \kern -.15em {\rm E}}}}
\newcommand{\bea}{\begin{eqnarray}}
\newcommand{\beas}{\begin{eqnarray*}}

\newcommand{\ba}{\begin{array}}
\newcommand{\ea}{\end{array}}

\newtheorem{theorem}{Theorem}[section]
\newtheorem{lemma}[theorem]{Lemma}

\newtheorem{remark}[theorem]{Remark}

\newcommand{\proof}{{\sc Proof.} \quad}
\newcommand{\eea}{\end{eqnarray}}                %
\newcommand{\bas}{\begin{eqnarray*}}            %
\newcommand{\eas}{\end{eqnarray*}}              %

 \journal{}
 
\begin{document}

% \begin{frontmatter}
%%%%%%%%%%%%%%%%%%%%%%%%%%%%%%%%%%%%%%%%%%%%%%%%%%%%%%%%%%%%%%%%%%%%%%%%%%%%%%%%%%%%%%%%%%%%%%%%%%%%%%%%%%%%%%%%%%%%%%%%%%%%

 \title{Title}

\author{Mihaela ~Negreanu}
\ead{negreanu@mat.ucm.es}
\address{Departamento de   Matem\'{a}tica Aplicada, Universidad Complutense de Madrid,  28040 Madrid, Spain}

\author{J. Ignacio ~Tello\corref{cor2}}
\ead{j.tello@upm.es}
\address{  Departamento   de  Matem\'{a}tica Aplicada a las tecnolog\'{\i}as  de la Informaci\'on  y las Telecomunicaciones,     E.T.S.I. Sistemas Inform\'aticos.
Universidad Polit\'ecnica de Madrid,  28031
Madrid, Spain \newline
Center for Computational Simulation, Universidad Polit\'ecnica de Madrid,
   28660 Boadilla del Monte, Madrid, Spain}
\tnotetext[ref1]{This work was supported by  Project MTM2013-42907 DGICT Spain}

\cortext[cor2]{Corresponding author}
% \fntext[fn1]{Telf.0034913945221}

\title{On a Parabolic-Elliptic   system with gradient dependent chemotactic coefficient}

 \date{}

\begin{abstract}
\noindent
We consider a second order PDEs  system of Parabolic-Elliptic type with chemotactic terms.  The system describes the evolution of a biological species ``$u$" moving
towards a higher concentration of a chemical stimuli  ``$v$" in a bounded and open  domain of $ \R^N$.  In the system considered,  the chemotaxis sensitivity depends on  the gradient of $v$, i.e., the chemotaxis term
has the following expression
$$- div  \left(\chi u |\nabla v|^{p-2}\nabla v \right),$$
where $\chi $ is a positive constant and $p$ satisfies
$$p \in (1, \infty), \quad \mbox{ if } N=1 \quad \mbox{ and } \quad p\in \left(1, \frac{N}{N-1}\right),   \quad \mbox{ if } N\geq 2.$$
We obtain uniform bounds in time in $L^{\infty}(\Omega)$ of the  solutions.
For the one-dimensional case we prove the existence of infinitely many non-constant steady-states for  $p\in (1,2)$ for any   $\chi$ positive and a given positive  mass.
% The paper concludes with some  remarks concerning the case  $p>2$.
\noindent
  \end{abstract}
\maketitle
{ \small  {\bf Keywords:}    Chemotaxis,  Global Existence of solutions, infinitely many solutions  }
%% keywords here, in the form: keyword \sep keyword

%% MSC codes here, in the form: \MSC code \sep code
%% or \MSC[2008] code \sep code (2000 is the default)

%\end{frontmatter}

%
%
%
%
%
%
%
%
%
%
%
%
%

\section{Introduction} \label{s0}
Chemotaxis is the   ability  of some  living organisms to orient their movement along a chemical concentration gradient. The process has been extensively  studied from a
biological point of view after the development of  the microscope during  the XIX century. In the last decades, several mathematical models have been presented to describe the phenomenon, after the pionnering works of Patlak \cite{P1}  and Keller and Segel \cite{KS}  (see also  the review articles  Horstmann \cite{H1},  \cite{H2}  and Bellomo et al \cite{bellomo1}
and references therein for more extensive literature in the subject).
The original model in \cite{KS} describes the evolution of a biological species, denoted by ``$u$" in terms of a parabolic equation,
with linear diffusion and a second order nonlinear term in the form
$$- div  ( \chi  u \nabla v),$$
where $v$ denotes the concentration of the chemical stimuli.
\\
In the last years, linear diffusion of the biological especies "$u$" has been replaced in different ways:
\begin{itemize}
\item by nonlinear diffusion at the form $-div ( \phi(u) \nabla u) $, see for instance Wrzosek \cite{wrzosek2010},  Cieslak and Mortales-Rodrigo \cite{cm},  Cieslak and Winkler \cite{cw}, Winkler \cite{winkler2017}.
\item by fractional diffusion,  see J Burczak, R Granero-Belinch\'on \cite{BG1} and \cite{BG2} among others.
\item by nonlinear diffusion depending on $|\nabla u|^{p}$,  (p-laplatian), see Bendahmane \cite{bendahmane}.
								\end{itemize}
The system has been also studied for several biological species, see   for instance
%  Lauffenburger \cite{l91},  Fasano et al \cite{fasano},  Conca and Espejo \cite{conca},
 Tang and Tao \cite{tt1}, Tello and Winkler \cite{tw1},   Stinner,  Tello  and   Winkler \cite{stw},
 Negreanu and Tello \cite{nt1} and \cite{nt2},  Wang  and Wu  \cite{ww} among others.
%\\
In the last years,  several mathematical models have considered the chemotactic sensitivity coefficient ``$\chi$"  dependent of $\nabla v$ instead of constant.
For instance, in Bellomo and Winkler \cite{bw1} and \cite{bw2} (see also Bellomo et al  \cite{bellomo2}), a chemotaxis  system is analyzed  for a chemotactic term of the form
 $$- div \left(  \frac{\chi  u}{[1+  |\nabla v|^2]^{\frac{1}{2}}}\nabla v\right).$$
In Bianchi, Painter and Sherratt \cite{bps}
the authors consider the term
$$-  div     \left(\frac{\chi u}{(1+\omega u)}  \frac{\nabla v }{(1+ \eta |\nabla v|)}\right) ,$$
for some positive constants $\chi$,  $\omega$ and $\eta$.
In \cite{bps},   the authors  study  a  system of four  PDEs of parabolic type  in an one-dimensional spatial domain  coupled with an ODE modeling Lymphangiogenesis in wound healing.
\\
A general chemotactic term is  presented  as
$$- div  \left[u \tilde{\chi}( u, v, |\nabla v|)\nabla v\right],$$
where $\tilde{\chi}$ is a continuos function for $|\nabla v|>0$.
\\
In the present work we consider a simplified  case, where $ \tilde{\chi}$ is given by
 $$ \tilde{\chi}( u, v, |\nabla v|)= \chi |\nabla v|^{p-2},$$
for   some positive constant $\chi $ and $p>1$.  % For the sake of  the notation we drop the tilde.
\\
We study the problem  in   a bounded spatial open domain $\Omega \subset \R^N$, with regular boundary $\partial \Omega$ and denote by $\vv{n}$ the  outward pointing normal vector on the boundary $\partial \Omega$.
%  and for simplicity   $$|\Omega|:= \int_{\Omega} dx =1$$.
The equation for $v$ is restricted to the elliptic case, for  simplicity, we assume that $v$ satisfies the Poisson equation and  the system studied is the following

\begin{equation}
\label{1}
 % \displaystyle  \left\{  \begin{array}{ll}   \displaystyle
 u_t-\Delta u= -  div  (\chi u|\nabla v|^{p-2}\nabla v), \quad \quad  x\in\Omega,\quad t>0,
     \end{equation}
     \begin{equation}
\label{1.1}
% \displaystyle
  -\Delta v  = u-M,\quad \quad \qquad \qquad \qquad \quad\;\; \;\;x\in\Omega,\quad t>0,
  %\end{array}  \right.
    \end{equation}
with   Neumann boundary conditions
\begin{equation}
 \displaystyle\frac{\partial  u}{\partial \vv{n} }=\displaystyle\frac{\partial v}{\partial \vv{n} }=0, \qquad \qquad \qquad\qquad \qquad \quad x\in \partial \Omega,\quad t>0  \label{2}
\end{equation}
and a non-negative    initial data
\begin{equation}
u(0,x)=u_0(x), \qquad \qquad \qquad \quad   x\in\Omega,\label{3}
\end{equation}
satisfying $$\frac{1}{|\Omega|}\displaystyle\int_{\Omega} u_0(x)dx=M.$$
Notice that if $(u,v)$ is a solution, $(u,v+k)$ is also a solution for any constant $k$.
To  obtain uniqueness  of the problem, we impose a given mass for the species $v$, i.e., we assume
\begin{equation}
\int_{\Omega} v=0.
\label{v0}
\end{equation}%
%  The problem (\ref{1}) - (\ref{3}) is equivalent to the following
%\begin{equation} \label{newq}  \left\{  \begin{array}{ll}  \displaystyle
% u_t-\Delta u= -  div  (\chi u \nabla w),  &   x\in\Omega,\quad t>0, \\[2mm]  \displaystyle
%  -div( |\nabla w|^{\frac{2-p}{p-1} } \nabla w ) = u-M, &    x\in\Omega,\quad t>0, \\[2mm]
% \displaystyle\frac{\partial  u}{\partial \vv{n} }=\displaystyle\frac{\partial w}{\partial \vv{n} }=0, &   x\in \partial \Omega,\quad t>0 \\[2mm]
%  u(0,x)=u_0(x), &  x\in\Omega,
%  \end{array} \right. \end{equation}
%  for $w$ defined by $$ \nabla w= |\nabla v|^{p-2}\nabla v.$$
%
The problem for $p=2$  has been already analyzed   by different authors, starting with  the work  of  J\"{a}ger and  Luckhaus \cite{JL} and   Biler
\cite{biler2}, and  \cite{biler1},
see also
Nagai \cite{nagai},  Senba \cite{senba},  Naito and Suzuki \cite{naito},   Blanchet, Dolbeault and  Perthame  \cite{blanchet2},  Blanchet, Carrillo and Masmoundi \cite{blanchet1} among others.
%
%Notice that the steady states of the system (\ref{newq}) satisfy
%\begin{equation} \label{newq2}  \left\{  \begin{array}{ll}  \displaystyle
%-\Delta u= -  div  (\chi u \nabla w),  &   x\in\Omega,  \\[2mm]  \displaystyle
%  -div( |\nabla w|^{\frac{2-p}{p-1} } \nabla w ) = u-M, &    x\in\Omega,  \\[2mm]
% \displaystyle\frac{\partial  u}{\partial \vv{n} }=\displaystyle\frac{\partial w}{\partial \vv{n} }=0, &   x\in \partial \Omega.
%\end{array} \right. \end{equation}
%Therefore $u= Me^{\chi w}$  and $w$ satisfies 
% \begin{equation} \label{newq3}  \left\{  \begin{array}{ll}  \displaystyle
 % -div( |\nabla w|^{\frac{2-p}{p-1} } \nabla w ) =M( e^{\chi w}-1), &    x\in\Omega, \\[2mm]
 % \displaystyle\frac{\partial w}{\partial \vv{n} }=0, &   x\in \partial \Omega,
%\end{array} \right. \end{equation}
%under the restriction
 %\begin{equation} \label{newq4}
%\frac{1}{|\Omega|} \int_{\Omega} e^{\chi w}=1.
%\end{equation}
%

In this article we analyze the   case  where $p$ satisfies
\begin{equation}\label{hipo00}p \in (1, \infty), \quad \mbox{ if } N=1 \quad \mbox{ and } \quad p\in \left(1, \frac{N}{N-1}\right),   \quad \mbox{ if } N\geq 2,  \end{equation}
%A priori estimates for the limit case $p=1$ can be obtained by using standard computations to deduce global existence of solutions. %A remark for the limit case $p=1$
% is included in section \ref{s2}
 assuming initial data
\begin{equation}
\label{a1} u_0 \in C^{2, \alpha} (\overline{\Omega}), \quad \mbox{ for some  $\alpha \in (0,1)$ }  \end{equation}
and
\begin{equation}
\label{a2}   \frac{\partial u_0}{\partial  \vv{n} } =0,  \quad x\in \partial \Omega.
\end{equation}
In Section \ref{s2} we study the global existence of the solutions $(u,v)$ of system (\ref{1})-(\ref{a2}). The main result is enclosed in   the following theorem.
\begin{theorem} Under assumptions (\ref{v0})-(\ref{a2}),   for any $T < \infty$,
there exists an unique classical solution  to (\ref{1})-(\ref{3}),
$$ u,v \in C^{2+\alpha, 1+\frac{\alpha}{2}}_{x,t}(\overline{\Omega}_T).$$
Moreover, there exists a constant $C(u_0, \chi, p, \Omega)$, independent of $T$,   such that
$$\|u\|_{L^{\infty}(\Omega)} \leq C .$$
\label{t1}
 \end{theorem}
%
% The proof of the theorem follows a Fixed Point Method  based in  a priori estimates  in $L^q(\Omega)$  of $u$.  % The  details are presented in Section \ref{s2}.
 \
 In Section \ref{s3} we consider  the steady states of the problem (\ref{1}) - (\ref{v0})  and we prove the existence of infinitely many solutions   in one dimensional bounded domain for $p\in (1,2)$. % The article finishes with a   discussion section where different open questions are treated.

\section{Global existence of Solutions} \label{s2}
\setcounter{equation}{0}
To prove  the global existence of solutions we apply
Schauder Fixed Point Theorem. We first introduce the local existence results  in Lemma \ref{lemma01}  and   obtain a priori estimates presented in the subsequent lemmas.  % In the following lemmas, we deduce some important properties of our solution.
\\

\begin{lemma}  \label{lemma01} Under assumptions (\ref{hipo00})-(\ref{a2}), there exists a unique  solution $(u,v)$ to
 (\ref{1})-(\ref{3})
in $(0,T_{max})$ satisfying
$$u, v \in C^{2+\alpha,  1+ \frac{\alpha}{2}}_{x,t} (\Omega_{T}), \qquad \mbox{ for any $T<T_{max}$}
$$
 where  $T_{max}$ is  a positive number satisfying 
 \begin{equation}
 \limsup_{t \rightarrow T_{max}} \left(  \|u(t)\|_{L^{\infty}(\Omega)}+ \|v(t)\|_{L^{\infty}(\Omega)}+ t  \right) =\infty.
    \label{b00}
\end{equation}
Moreover, the solution $u$  satisfies
\begin{equation} \label{positive}
u(t,x)\geq 0, \qquad x\in\Omega, \quad t<T_{max}. \end{equation}
\end{lemma}

{\bf Proof.} For any $T<T_{max}$ we have that $u$ satisfies
$$u_t-\Delta u +b(t,x)\cdot  \nabla u =  f(x,t), \quad (t,x) \in \Omega_T,$$ 
where $$b(x,t)=  \chi  |\nabla v|^{p-2}\nabla v, \qquad f(x,t)=div  \left(\chi  |\nabla v|^{p-2}\nabla v \right).$$
Since $u \in L^{\infty}(\Omega_T)$ we have that $v \in L^s(0,T: W^{2,q}(\Omega)) $ for any $s, q\geq 1$ and therefore 
 $$b(x,t)\in L^s(0,T: W^{1,q}(\Omega)), \qquad f(x,t)\in L^s(0,T: L^{q}(\Omega)).$$
Then, $u \in C^{2+\alpha,  1+ \frac{\alpha}{2}}_{x,t} (\Omega_{T})$ see Remark 48.3 (ii) in Quittner-Souplet \cite{qs}.

% The proof of local existence follows standard semigroup theory and fixed point argument, we refer to   (see also   Horstmann and Winkler \cite{hw1},  Tello and  Winkler \cite{tw01} and Tao and Winkler \cite{tao}  for more details).
 The non-negativity of $u$ is a consequence of the maximum principle.   \qed
%and maximum principle for elliptic equation  end the proof.     % The solution is extended to a maximum interval of existence $(0, T_{max})$.

Notice that, after integration in (\ref{1}), we have that the total mass is preserved in time, i.e.,
\begin{equation} \label{totalmass}
\int_{\Omega} u= \int_{\Omega} u_0=|\Omega| M. 
\end{equation}

 \begin{remark} \label{R01}  Let $\Omega \subset \R^N$, be a bounded and regular domain and $f\in L^1(\Omega)$, such that $\int_{\Omega} f=0$,  then, the problem
$$ \left\{
  \begin{array}{l}
\displaystyle\displaystyle- \Delta v= f, \qquad \mbox{ in }  \Omega
 \\ \\
\displaystyle \frac{\partial v}{\partial  \vv{n}}=0, \quad \quad\;\;\mbox{ in }   \partial \Omega,
 \end{array}   \right.
 $$
 has an unique solution
 \begin{equation} v \in W^{1, s}(\Omega), \quad \mbox{ for } \quad  \left\{ \begin{array}{l}   s\in [1, N/(N-1)),  \mbox{ if  }  N>1 \\[2mm] s= \infty  \mbox{ if }  N>1 ,  \end{array}  \right. \label{lemmaPE}  \end{equation}
  satisfying
$$ \int_{\Omega } v=0.$$
The proof of (\ref{lemmaPE}) is given in Chabroswki \cite{chabrowski} Theorem 2.8, where the general problem
$$ \left\{
  \begin{array}{l}
\displaystyle- \Delta v= \lambda v + f, \qquad \mbox{ in }  \Omega
 \\ \\
\displaystyle \frac{\partial v}{\partial   \vv{n} }=0, \quad\qquad \qquad\; \mbox{ in }   \partial \Omega,
 \end{array}   \right.
 $$
 is  studied for $\lambda  \in \R$.
  \end{remark}

\begin{lemma} Let $p>1$ such that
 \begin{equation}
  \label{hq1} \left\{ \begin{array}{l} p<\infty, \quad\quad \quad \quad \quad  \mbox{ if } \quad N=1,
\\[3mm] 2(p-1)< \frac{N}{N-1},    \quad \mbox{ if } \quad  N\geq 2,    \end{array} \right. \end{equation}
 then, for  any
 $q>1$   and any  $s>0$  satisfying   \begin{equation}
  \label{hq1b} \left\{ \begin{array}{c}  s\in \left(2(p-1),  \infty \right), \quad \quad \;\mbox{ if } \quad N=1,  \\[2mm]
  s\in \left(2(p-1),  \frac{N}{N-1}\right), \quad \mbox{ if } \quad N\geq 2,  \end{array} \right. \end{equation}
   the following inequality holds:
      \begin{equation}
\label{5.1}
 \frac{d}{dt}\int_{\Omega}u^q+\frac{3(q-1)}{q} \int_{\Omega}\left|\nabla u^{\frac{q}{2}}\right|^2 \leq c_1 q (q-1) \chi^2
  \left[  \int_{\Omega} u^{\frac{sq}{s-2(p-1)}} \right]^{\frac{s-2(p-1)}{s}} ,
    \end{equation}
    for some positive constant $c_1$.
 \label{l0}
\end{lemma}
{\bf Proof:}
 We multiply equation (\ref{1}) by $u^{q-1}$  (for $q>1$) and integrate by parts to obtain
\begin{equation}
\label{4}
  \displaystyle  \frac{d}{dt}\frac{1}{q}\int_{\Omega}u^q+(q-1)\int_{\Omega}\left|\nabla u\right|^2 u^{q-2}= (q-1) \chi\int_{\Omega}\left|\nabla v\right|^{p-2}\nabla v \nabla u \cdot u^{q-1}.
    \end{equation}
    Thanks to Young's inequality, it results
   $$
     \begin{array}{ll}
&\displaystyle \frac{d}{dt}\frac{1}{q}\int_{\Omega}u^q+\frac{4(q-1)}{q^2}\int_{\Omega}\left|\nabla u^{\frac{q}{2}}\right|^2\displaystyle \leq \frac{2 (q-1)\chi }{q}\int_{\Omega} \left|\nabla v\right|^{p-1} \left|\nabla u^{\frac{q}{2}}\right|  u^{\frac{q}{2}}\\
\noalign{\medskip}\displaystyle &
\leq \displaystyle\frac{ (q-1)}{q^2} \int_{\Omega} \left|\nabla u^{\frac{q}{2}}\right|^2+(q-1)\chi^2
\int_{\Omega} \left|\nabla v \right|^{2p-2} u^{q}
\\
\noalign{\medskip  }\displaystyle &\displaystyle\leq
\frac{(q-1)}{     q^2}\int_{\Omega} \left|\nabla u^{\frac{q}{2}}\right|^2 + (q-1) \chi^2
\left[\int_{\Omega} \left|\nabla v\right|^{s}\right]^{\frac{2(p-1)}{s}} \left[  \int_{\Omega} u^{\frac{sq}{s-2(p-1)}} \right]^{\frac{s-2(p-1)}{s}},
 \end{array}
 $$
  then \begin{equation} \label{qsq}
  \frac{d}{dt}\frac{1}{q}\int_{\Omega}u^q+\frac{3(q-1)}{q^2}\int_{\Omega}\left|\nabla u^{\frac{q}{2}}\right|^2\displaystyle \leq(q-1) \chi^2
\left[\int_{\Omega} \left|\nabla v\right|^{s}\right]^{\frac{2(p-1)}{s}} \left[  \int_{\Omega} u^{\frac{sq}{s-2(p-1)}} \right]^{\frac{s-2(p-1)}{s}}.
\end{equation}
We notice that, since $u \in  L^1(\Omega)$, we have $v\in W^{1,s}(\Omega)$ for any $s$ satisfying (\ref{hq1b})
 and, thanks to Lemma  \ref{lemmaPE}, we have
\begin{equation} \|\nabla v\|_{L^s(\Omega)}^{2(p-1)} \leq c_1<\infty. \label{bre} \end{equation}
Then,
equation (\ref{qsq}) becomes
$$
 \frac{d}{dt}\int_{\Omega}u^q+ \frac{3(q-1)}{q}\int_{\Omega}\left|\nabla u^{\frac{q}{2}}\right|^2 \leq  \chi^2(q-1)q c_1
  \left[  \int_{\Omega} u^{\frac{sq}{s-2(p-1)}} \right]^{\frac{s-2(p-1)}{s}},
$$
 the  proof is done.    \qed

    \begin{lemma}
  \label{l4} Under assumption \begin{equation}
  \label{hq2} \left\{ \begin{array}{l} p\in (1, \infty ),\qquad \quad \;\;\quad \mbox{ if } N=1,
\\[3mm]
p \in \left(1, \frac{N}{N-1} \right), \qquad \quad \mbox{ if } N\geq 2,    \end{array} \right. \end{equation}
   and any $s>0$ satisfying
  \begin{equation}
  \label{hq2b} \left\{ \begin{array}{c}  s\in \left(2(p-1),  \infty \right), \quad \qquad \mbox{ if } \quad N=1,  \\[2mm]
  s\in \left(2(p-1),  \frac{N}{N-1}\right), \quad \;\mbox{ if } \quad N\geq 2,  \end{array} \right. \end{equation}
 we have that \begin{equation}
    \label{12bis}\int_{\Omega} u^{N+1} \leq c_5, \qquad \mbox{ for any } t>0,\end{equation}
and  $c_5$  independent of $t$.
  \end{lemma}
 \proof We first notice that for $N \geq 2$     we have that $p<2$    and then  $p \geq 2(p-1)$ and, in view of (\ref{hq2}),
 $$2(p-1) \leq p<\frac{N}{N-1}$$
 which gives a non empty set  of admissible  values of  $s$ for $N\geq 2$.
\\
 To prove the Lemma we follow an Moser-Alikakos iteration method, the result could be also  obtained using similar arguments as in Tao and Winkler \cite{tao}, for readers convenience we detail the proof.   We recall the expression    (\ref{5.1})
\begin{equation*}
  \frac{d}{dt}\int_{\Omega}u^q+ \frac{3(q-1)}{q}\int_{\Omega}\left|\nabla u^{\frac{q}{2}}\right|^2 \leq  \chi^2(q-1)q c_1
  \left[  \int_{\Omega} u^{\frac{sq}{s-2(p-1)}} \right]^{\frac{s-2(p-1)}{s}}.
\label{qqq} \end{equation*}
%
%  If  $2(p-1)<\frac{N}{N-1}$  then
%
By Gagliardo-Nirenberg  inequality  (see Henry \cite{Henry}),  for $1\leq \gamma \leq  \beta \leq \infty $,  $ r\geq 1$
and  $$-\frac{N}{\beta}=  a (1-\frac{N}{r})-\frac{N}{\gamma}(1- a),$$ if  $$a\in (0,1),$$
 we have that
\begin{equation} \label{gagli} \|w \|_{L^{\beta}(\Omega)} \leq C_{GN} \|w \|_{L^{\gamma}(\Omega)}^{1-a} \|  w  \|_{W^{1,r}(\Omega)   }^{a} . \end{equation}
Notice that
$$a=\frac{\frac{1}{\gamma}- \frac{1}{\beta} }{\frac{1}{\gamma}+\frac{1}{N}- \frac{1}{r}}.
$$
We take $$w=u^{\frac{q}{2}}, \quad \gamma=1  \qquad  \beta= \frac{2s}{   s-2(p-1)}, \quad r=2, \quad \mbox{ and } \  s= \frac{N}{N-1+\epsilon}$$
for  $\epsilon $ small enough  such that  \begin{equation} s\in \left(2(p-1),  \frac{N}{N-1}\right), \quad 
N- \frac{p}{p-1}+ \epsilon <0 .\label{ep12} \end{equation}   Then, we  obtain
$$\|u^{\frac{q}{2}    }\|_{L^{\frac{2s}{   s-2(p-1)} }(\Omega)} \leq C_{GN} \| u^{\frac{q}{2} } \|_{H^1(\Omega)}^{a} \|u^{\frac{q}{2} }\|_{L^1(\Omega)}^{1-a},  $$
for   $$\begin{array}{ll}
\displaystyle a&= \displaystyle\frac{1-  \frac{ s-2(p-1)  }{2s}   }{ \frac{1}{2}+ \frac{1}{N}}
 =\frac{2N}{N+2}  (1-  \frac{  s-2(p-1)  }{2s}   ) = \frac{N}{N+2}  \frac{ (s+2(p-1) ) }{s}   \\
\noalign{\medskip}\displaystyle &\displaystyle= \frac{N}{N+2}  \frac{ (\frac{N}{N-1+\epsilon}+2(p-1) ) }{\frac{N}{N-1+\epsilon}},   \end{array}
 $$
 equivalent to
 $$\begin{array}{ll}
\displaystyle a&= \displaystyle\frac{N-1 + \epsilon }{N+2}   \left(\frac{N}{N-1+\epsilon}+2(p-1) \right)  =    \frac{N + 2(p-1)(N-1+\epsilon)}{N+2}\\ \noalign{\medskip}\displaystyle &\displaystyle<    \frac{N +2 + 2(p-1)N-2p+2(p-1)\epsilon}{N+2}
 = 1 +2(p-1)\frac{ N-\frac{p}{p-1}+\epsilon}{N+2} .  \end{array}   $$
Notice that $$N< \frac{p}{p-1}\quad \quad\quad\mbox{ $\Longleftrightarrow$ } \quad\quad\quad
 p< \frac{N}{N-1},$$
and in view of  (\ref{ep12})  we have $a<1$.

 We now apply inequality (\ref{gagli})  to the last term in (\ref{5.1}) to obtain
 $$  \left[  \int_{\Omega} u^{\frac{sq}{s-2(p-1)}} \right]^{\frac{s-2(p-1)}{s}} \leq  C_{GN} \left[ \| u^{\frac{q}{2} } \|_{H^1(\Omega)}^{a} \|u^{\frac{q}{2} }\|_{L^1(\Omega)}^{1-a} \right]^{2}
 $$
and
 $$C_{GN} \left[ \| u^{\frac{q}{2} } \|_{H^1(\Omega)}^{a} \|u^{\frac{q}{2} }\|_{L^1(\Omega)}^{1-a} \right]^{\frac{s-2(p-1)}{s}} \leq \delta \| u^{\frac{q}{2} } \|_{H^1(\Omega)}^{2}  + \frac{c_2}{\delta^{\frac{a}{(1-a)}}}
 \|u^{\frac{q}{2} }\|_{L^1(\Omega)}^{2}  ,
 $$
for some positive $\delta $ small enough. Using the above estimates, (\ref{5.1}) becomes
$$
  \frac{d}{dt}\int_{\Omega}u^q+ \frac{3(q-1)}{q}\int_{\Omega}\left|\nabla u^{\frac{q}{2}}\right|^2 \leq  \chi^2(q-1)q c_1
  \left[  \delta  \int_{\Omega}\left|\nabla u^{\frac{q}{2}}\right|^2 + \delta  \int_{\Omega}u^{q}+  \frac{c_2}{\delta^{\frac{a}{(1-a)}}} \|u^{\frac{q}{2} }\|_{L^1(\Omega)}^{2}   \right],
$$
 i.e.,
$$
  \frac{d}{dt}\int_{\Omega}u^q+ \left[\frac{3(q-1)}{q}- \delta   \chi^2(q-1)q c \right] \int_{\Omega}\left|\nabla u^{\frac{q}{2}}\right|^2 \leq  \chi^2(q-1)q c
  \left[  \delta  \int_{\Omega}u^{q}+ \frac{c_2}{\delta^{\frac{a}{(1-a)}}} \|u^{\frac{q}{2} }\|_{L^1(\Omega)}^{2}   \right]
$$
and for $\delta < \frac{2}{\chi^2 c q^2}$ it results
 \begin{equation}
  \frac{d}{dt}\int_{\Omega}u^q+  \frac{(q-1)}{q}  \int_{\Omega}\left|\nabla u^{\frac{q}{2}}\right|^2 \leq  \chi^2(q-1)q c
  \left[  \delta  \int_{\Omega}u^{q}+ \frac{c_2}{\delta^{\frac{a}{(1-a)}}}  \|u^{\frac{q}{2} }\|_{L^1(\Omega)}^{2}   \right].
\label{qqq2} \end{equation}
Thanks to  Poincar\'e-Wirtinger inequality we get 
 $$  \int_{\Omega}\left|\nabla u^{\frac{q}{2}} \right|^2  \geq C_{PW}    \int_{\Omega}  \left(   u^{\frac{q}{2}}  - \frac{1}{|\Omega| } \int_{\Omega}  u^{\frac{q}{2}}  \right)^2 =
 C_{PW}   \left(   \int_{\Omega}   u^{q}  - \frac{1}{|\Omega| } \left|\int_{\Omega}  u^{\frac{q}{2}}  \right|^2 \right),
  $$
which implies
$$
\frac{d}{dt} \int_{\Omega}u^q+\frac{q-1}{q}\left[\frac{  C_{PW} }{2}-\delta q^2 c \chi^2  \right] \int_{\Omega}  u^{q}    \leq   c_3(\delta, |\Omega| , q) \left|\int_{\Omega}  u^{\frac{q}{2}}  \right|^2   .  $$
For $\delta <  \min\{ \frac{C_{PW}}{4c \chi^2 q^2},   \frac{2}{\chi^2 c q^2} \} $ the last inequality is reduced to
\begin{equation} \label{33e}
\frac{d}{dt} \int_{\Omega}u^q+\frac{q-1}{q}\frac{  C_{PW} }{4}  \int_{\Omega}  u^{q}    \leq   c_3(\delta, |\Omega| , q) \left|\int_{\Omega}  u^{\frac{q}{2}}  \right|^2   \end{equation}
and by the Maximum Principle,  we get
\begin{equation} \label{33f} \sup_{t>0} \|u\|_{L^q(\Omega)}^q \leq \frac{4q c_3(\delta,  |\Omega|, q) }{(q-1) C_{PW}}  \sup_{t>0}  \|u\|^{q}_{L^{\frac{q}{2}}(\Omega)} . \end{equation}
Following Moser-Alikakos iteration (see \cite{alikakos}), we define
$$x_{i}:= \sup_{t >0}   \int_{\Omega} u^{2^i}, \qquad i \in \N.$$
Thanks to
(\ref{totalmass})  it results 
$x_0=|\Omega| M<\infty$  and (\ref{33f}) implies that
$$x_i\leq c_4(\delta,i) x_{i-1}^2.$$
Notice that    $x_i$ is finite for $i<\infty$, in particular, there exists $i_0$ large enough, such that $2^{i_0}>N+1$  and  therefore
$$ \int_{\Omega} u^{N+1}\leq c_5<\infty \quad \mbox{ for any } t>0,$$ and the proof ends.
\qed
\begin{lemma}  \label{lqw} There exists a positive constant $c_{6}$, independent of $t$, such that the following bound holds
$$\|\nabla v\|_{L^{\infty}(\Omega)} <c_6.$$
\end{lemma}
{\bf Proof:} In view of (\ref{12bis}) we claim that $v\in W^{2,N+1}(\Omega) \subset W^{1, \infty}(\Omega)$, thanks to  
 Chabroswki \cite{chabrowski} Theorem 2.8,  we have the result.
\qed
  \begin{lemma} \label{l6} There exists a positive constant $c_{\infty}$, independent of $t$, such that
  $$ \|u\|_{L^{\infty}(\Omega) }\leq c_{\infty}
  .$$
\end{lemma}
{\bf Proof:}
From
(\ref{4}) (for $q\geq 2$) and Lemma \ref{lqw} we have
$$
 \frac{d}{dt}\frac{1}{q}\int_{\Omega}u^q+(q-1)\int_{\Omega}\left|\nabla u\right|^2 u^{q-2}= (q-1) c_7\int_{\Omega}   |\nabla u |  u^{q-1}.
  $$
Since
$$
c_7 \int_{\Omega}   |\nabla u |  u^{q-1} \leq  \frac{1}{2} \int_{\Omega}\left|\nabla u\right|^2 u^{q-2}+  \frac{c_7^2}{2} \int_{\Omega} u^{q},
$$
we get
\begin{equation}
 \frac{d}{dt}\frac{1}{q}\int_{\Omega}u^q+\frac{(q-1)}{2}\int_{\Omega}\left|\nabla u\right|^2 u^{q-2}= (q-1) \frac{c_7^2}{2}\int_{\Omega}    u^{q}.
  \label{eq44}
  \end{equation}
In the Gagliardo-Nirenberg's inequality
$$\|w\|_{L^2(\Omega)}\leq  C_{GN1} \|\nabla w\|^a_{L^2(\Omega)} \| w\|^{1-a}_{L^1(\Omega)} + C_{GN2}\|w\|_{L^1(\Omega)}, \quad \mbox{ for }  \quad
\frac{1}{2}= a\left( \frac{1}{2}- \frac{1}{N}\right) +1-a,
$$
taking $$w= u^{\frac{q}{2}}, \qquad \mbox{ for } \quad  a=\frac{N}{N+2}<1, $$  we obtain
   $$\left[ \int_{\Omega}    u^{q} \right]^{\frac{1}{2}} \leq  C_{GN1}  \left[\int_{\Omega} |\nabla u^{\frac{q}{2}} |^2\right]^{\frac{a}{2}} \left[ \int_{\Omega}    u^{\frac{q}{2}} \right]^{(1-a)} +  C_{GN2} \int_{\Omega}    u^{\frac{q}{2}}, $$
 which is  equivalent to
    $$ \int_{\Omega}    u^{q}  \leq 2 C_{GN1}   \left[\int_{\Omega} |\nabla u^{\frac{q}{2}} |^2\right]^{a} \left[ \int_{\Omega}    u^{\frac{q}{2}} \right]^{2(1-a)} +2  C_{GN2}  \left[ \int_{\Omega}    u^{\frac{q}{2}} \right]^{2}.$$
    Replacing the last expression in (\ref{eq44}) we have
\begin{equation} \label{beq44}
 \frac{d}{dt}\frac{1}{q}\int_{\Omega}u^q+\frac{2(q-1)}{q^2}\int_{\Omega}\left|\nabla u^{\frac{q}{2}}\right|^2 \leq  (q-1) c_8\left(
\left[\int_{\Omega} |\nabla u^{\frac{q}{2}} |^2\right]^{a} \left[ \int_{\Omega}    u^{\frac{q}{2}} \right]^{2(1-a)} +\left[ \int_{\Omega}    u^{\frac{q}{2}} \right]^{2} \right)
 \end{equation}
   and applying Young's inequality to the first term of the right side, 
$$\left[\int_{\Omega} |  \nabla u^{\frac{q}{2}} |^2\right]^{a} \left[ \int_{\Omega}    u^{\frac{q}{2}} \right]^{2(1-a)}     \leq
\frac{1}{q^2c_8} \int_{\Omega} |\nabla u^{\frac{q}{2}} |^2 + \left[q^2c_8a\right]^{\frac{a}{(1-a)}} (1-a)   \left[ \int_{\Omega}    u^{\frac{q}{2}} \right]^{2}, $$
  inequality (\ref{beq44}) is reduced to
   $$\frac{1}{q}
 \frac{d}{dt}\int_{\Omega}u^q+\frac{(q-1)}{q^2}\int_{\Omega}\left|\nabla u^{\frac{q}{2}}\right|^2 \leq  (q-1)
 q^{\frac{2a}{(1-a)}}  c_9  \left[ \int_{\Omega}    u^{\frac{q}{2}} \right]^{2},
$$
    i.e.,
     $$\frac{q}{(q-1)}
 \frac{d}{dt}\int_{\Omega}u^q+ \int_{\Omega}\left|\nabla u^{\frac{q}{2}}\right|^2 \leq
 q^{\frac{2}{(1-a)}}  c_9  \left[ \int_{\Omega}    u^{\frac{q}{2}} \right]^{2}
$$
Thanks to  Poincar\'e-Wirtinger's inequality we have
 $$  \int_{\Omega}\left|\nabla u^{\frac{q}{2}} \right|^2  \geq C_{PW}    \int_{\Omega}  \left(   u^{\frac{q}{2}}  - \frac{1}{|\Omega| } \int_{\Omega}  u^{\frac{q}{2}}  \right)^2 =
 C_{PW}   \left(   \int_{\Omega}   u^{q}  - \frac{1}{|\Omega| } \left|\int_{\Omega}  u^{\frac{q}{2}}  \right|^2 \right),
  $$
which implies
$$
\frac{q}{q-1} \frac{d}{dt} \int_{\Omega}u^q+ C_{PW}   \int_{\Omega}  u^{q}    \leq   q^{2+\frac{2a}{(1-a)}}  c_{10}   \left|\int_{\Omega}  u^{\frac{q}{2}}  \right|^2    .  $$
By the Maximum Principle  we get
\begin{equation} \label{33qer} \sup_{t>0} \|u\|_{L^q(\Omega)}^q \leq \frac{   q^{\frac{2}{(1-a)}}  c_{10}   }{ C_{PW}}   \sup_{t>0}  \|u\|^{q}_{L^{\frac{q}{2}}(\Omega)} . \end{equation}
As in Lemma \ref{l4} we apply the following Moser-Alikakos iteration (see \cite{alikakos}), defining
$$x_{i}:= \sup_{t >0}   \int_{\Omega} u^{2^i}, \qquad i \in \N.$$
Thanks to
(\ref{totalmass})  we have that
$x_0=|\Omega| M<\infty$  and (\ref{33f}) implies that
$$x_i\leq   \displaystyle 2^{\frac{2 i }{(1-a)}}  c_{11} x_{i-1}^2,$$
then
$$x_i^{2^{-i}} \leq   \displaystyle 2^{\frac{i2^{1-i}}{(1-a)}  }  c_{11}^{2^{-i}}  x_{i-1}^{2^{1-i}},$$
and
$$\displaystyle\|u\|_{L^{2^i}(\Omega)} \leq   \displaystyle\prod_{j=0\dots i}   2^{\displaystyle \frac{j2^{1-j}}{(1-a)}   }  \displaystyle c_{11}^{2^{-j}} = \displaystyle2^{\displaystyle\sum_{j=0 \dots i} \displaystyle\frac{j2^{1-j}}{(1-a)} }  c_{11}^{\displaystyle\sum_{j=0 \dots i}  2^{-j} }.  $$
Since
$$\begin{array}{ll}
\displaystyle \sum_{j=0 \dots i} \frac{j 2^{1-j}}{(1-a)} &= \displaystyle   \frac{2}{(1-a)}  \sum_{j=0 \dots i} j2^{-j}= \frac{2}{(1-a)} 
(2-2^{-i}-2^{-i-1}) \leq  \frac{4 }{(1-a)}  \end{array} $$
and 
$$ \sum_{j=0 \dots i}  2^{-j}  = 2- 2^{-i}<2,$$
we obtain
$$ \|u\|_{L^{2^i}(\Omega)} \leq    2^{\frac{2}{(1-a)}} c_{11}^2:= c_{\infty}.$$
Constant  $c_{11}$ is independent of $i$,
taking limits when $i \rightarrow \infty$, the proof ends. \qed
\\
{\bf End of the proof of Theorem \ref{t1}.}
\newline
In view  of Lemma \ref{lqw} and assumption \ref{v0} we obtain uniform boundedness of $v$ in $L^{\infty}(\Omega_{T_{max}})$. Thanks to Lemmatta  \ref{lemma01}  and \ref{l6} the proof ends.

\begin{remark}
We notice that, if equation (\ref{1.1}) is replaced by
   \begin{equation}
\label{1.1R}
% \displaystyle
  -\Delta v +v = u,\quad \quad  x\in\Omega,\quad t>0,
  %\end{array}  \right.
    \end{equation}
and Remark \ref{R01}  is replaced by  Lemma 23 in  Brezis-Strauss \cite{brezis}, we  also obtain  the boundedness in $L^{\infty}(\Omega)$ of $u$ and $v$ and the global existence of solutions.
\end{remark}

\section{Stationary states in 1 D} \label{s3}
\setcounter{equation}{0}
The steady states of the problem (\ref{1})-(\ref{v0})  for $\Omega=(0,1)$, are given by

\begin{equation} \label{3.1}\left\{ \begin{array}{ll}  - u_{xx}= -   (\chi u|v_x|^{p-2}v_x)_x , &    x\in(0,1),  \\
    -v_{xx}  = u-M,&   x\in(0,1),   \\
  u_x(0)=u_x(1)=v_x(0)=v_x(1)=0.  &  \end{array} \right. \end{equation}
The aim of this section is to prove the following theorem:
\begin{theorem} Let $p\in (1,2)$, then,
for any $M>0$ and $\chi>0$, there exist infinitely many  solutions to (\ref{3.1}).
\label{T3.1}
\end{theorem}
In order to prove the theorem, we proceed into several steps.
\\
We consider the following  variables   \begin{equation}
\label{18}w: = ln(u)-ln(M), \qquad \beta: =v_x \end{equation}
then, (\ref{3.1})  becomes
\begin{equation}
\label{19}
 \displaystyle
  \left\{
  \begin{array}{ll}
  \displaystyle  w_x= \chi |\beta|^{p-2}\beta, \quad \quad & x\in(0,1),
     \\
\displaystyle   \beta_x=M( 1-e^w),\quad \quad & x\in (0,1),  \\  \beta(0)=\beta(1)=0.  &  \end{array}
  \right.
    \end{equation}
  We have that (\ref{19}) is a Hamiltonian system, i.e., there exists $H(w,\beta)$ such that
\begin{equation}
\label{19bis}
 \displaystyle
  \left\{
  \begin{array}{ll}
  \displaystyle  \chi |\beta|^{p-2}\beta=\frac{\partial H(w,\beta)}{\partial\beta},\quad \quad & x\in (0,1)
     \\
\displaystyle   M( 1-e^w)=-\frac{\partial H(w,\beta)}{\partial w},\quad \quad & x\in (0,1).  \end{array}
  \right.
    \end{equation}
After integration,   we obtain
\begin{equation}
\label{20} H(w,\beta)=M( e^w-w-1)+  \chi\frac{ |\beta|^{p}}{p}. \end{equation}
The solutions  of the system are found along the contours of $H$, so that, to sketch the phase diagram of this system, it is enough to study and draw the level sets of the Hamiltonian function $H$. We notice that  there exists an energy in the  system (\ref{19}) which is preserved along the solutions, i.e.,
\begin{equation}
\label{21} \frac{M}{\chi}(e^w-w-1)+ \frac{1}{p}|\beta|^p=k,\end{equation}
where $k$ is a non negative constant.  \ 
Notice that  for   $k= 0$ we have   the trivial solution   $w=\beta=0$.
\\
In the following lemma we prove the conservation of the energy.
\begin{lemma}
\label{lema21}
Let $(w, \beta)$ be a solution  to (\ref{19}) and the hamiltonian function $H$ as in (\ref{20}), i.e.,
\begin{equation*}
\label{22}H(w, \beta)=M(e^w-w-1)+\frac{\chi}{p}|\beta|^p \end{equation*}
 then,
$$\frac{d H}{dx} =0.$$
\end{lemma}
{\bf Proof:} Recalling the Hamilton function $H$ is conserved in any solution of the system, i.e.,
$$\frac{d H(w(x),\beta(x))}{dx}=\frac{\partial H}{\partial w}\frac{\partial w}{\partial x}+\frac{\partial H}{\partial \beta}\frac{\partial \beta}{\partial x}=\frac{\partial H}{\partial w}\frac{\partial H}{\partial \beta}-\frac{\partial H}{\partial \beta}\frac{\partial H}{\partial w}=0,$$
and we have the proof. \qed
\\
Notice that the level sets of $H$ are bounded curves corresponding to a periodic solution $(w, \beta)$ of period $T(k)$,  for $k$ defined  in  (\ref{21}).

{\bf End of proof of Theorem  \ref{T3.1}.}
\\
%  We introduce the function  $\phi: (0,\infty)\rightarrow (0, \infty)$,  defined by
%  \begin{equation} \label{phi} \phi(s)=e^{-s}+s-1 \end{equation}
%  which  is  a monotone increasing function for $s >0$. % Function $\phi$ will be used in the following lemma.
We write  $|\beta|$ in terms of $|w_x|$, then
$$|\beta|^{p}= \chi^{\frac{p}{1-p}} |w_x|^{\frac{p}{p-1}}$$
and
\begin{equation}
\label{21c} \frac{M}{\chi}(e^w-w-1)+ \frac{\chi^{\frac{p}{1-p}}}{p}|w_x|^{\frac{p}{p-1}}=k\end{equation}
 and $|w_x|$ satisfies 
 \begin{equation} \label{22}|w_x|=  [kp \chi^{\frac{p}{p-1}}+ Mp \chi^{\frac{1}{p-1}}(1+w-e^w)]^{\frac{p-1}{p}}.\end{equation}
We denote by $r_0$ and $-r_1$ the values of $w$ for $\beta =0$ for a given $k$, then  $r_1$ and $r_0$ satisfy
$$ k= \frac{M}{\chi} (e^{r_0}-r_0-1)=\frac{M}{\chi} (e^{-r_1}+r_1-1). $$
\\
After integration in (\ref{22})  we get
$$ \int_{-r_1}^{r_0} \left[kp \chi^{\frac{p}{p-1}}+ Mp \chi^{\frac{1}{p-1}}(1+w-e^w)\right]^{\frac{1-p}{p}} dw= \frac{T}{2},$$
i.e. 
\begin{equation} \label{int} \int_{-r_1}^{r_0} \left[ \frac{k \chi}{M} +(1+w-e^w)\right]^{\frac{1-p}{p}} dw=[Mp\chi^{\frac{1}{p-1}}]^{\frac{p-1}{p}} \frac{T}{2}, \end{equation}
and 
$$\int_{-r_1}^{r_0} \left[\frac{k \chi}{M} +(1+w-e^w)\right]^{\frac{p-1}{p}} dw =
  \int_{-r_1}^{0}  \left[ e^{-r_1}+r_1-e^w+w \right]^{-\frac{p-1}{p}}  dw$$
$$\qquad\qquad\qquad\qquad\qquad\qquad\qquad\qquad+
  \int_{0}^{r_0} \left[ e^{r_0}-r_0 -e^w+w\right]^{-\frac{p-1}{p}} dw.$$
To estimate the righthand side integrals, we apply H\"older inequality, so, for any $\epsilon>0$  verifying
\begin{equation} \label{epsilon}  \epsilon <\frac{2-p}{p-1} \end{equation}  we have:
 $$  \int_{0}^{r_0} \left[ e^{r_0}-r_0 -e^w+w \right]^{-\frac{p-1}{p}} dw\leq  \int_{0}^{r_0} (e^{w}-1)  \left[ e^{r_0}-r_0 -e^w+w  \right]^{-\frac{(p-1)(2+ \epsilon)}{p}  } dw
   +  \int_{0}^{r_0} (e^{w}-1)^{-\frac{1}{1+ \epsilon}}.
   $$
   After integration, for $p<2$, we get
    $$  \int_{0}^{r_0}(e^{w}-1) \left[ e^{r_0}-r_0 -e^w+w\right]^{-\frac{(p-1)(2+\epsilon)}{p}} dw=  \frac{p}{(2+ \epsilon)-p(1+ \epsilon)}\left[ e^{r_0}-r_0-1 \right]^{\frac{(2+ \epsilon)-p(1+ \epsilon)}{p}  } $$
    and
    $$ \int_{0}^{r_0} (e^{w}-1)^{-\frac{1}{1+ \epsilon}} \leq  \int_{0}^{r_0} w^{-\frac{1}{1+ \epsilon} } = \frac{1+\epsilon}{ \epsilon}
    r_0^{ \frac{\epsilon}{1+ \epsilon}} .$$
    In the same way  we obtain
    $$ \int_{-r_1}^{0}  \left[ e^{-r_1}+r_1-e^w+w \right]^{-\frac{p-1}{p}}  dw \leq \int_{-r_1}^{0} (1-e^w)  \left[ e^{-r_1}+r_1-e^w+w \right]^{-\frac{(p-1)(2+\epsilon)}{p}}  dw $$
   $$\qquad\qquad+ \int_{-r_1}^{0}     (1-e^{w})^{-\frac{1}{1+ \epsilon}}  dw.$$
 Computing the above integrals it follows
   $$\int_{-r_1}^{0} (1-e^w)  \left[ e^{-r_1}+r_1-e^w+w \right]^{-\frac{(p-1)(2+\epsilon)}{p}}  dw =
   \frac{p}{(2+ \epsilon)-p(1+ \epsilon)}  \left[ e^{-r_1}+r_1-1 \right]^{1-\frac{(p-1)(2+\epsilon)}{p}}
   $$ and
     $$ \int_{-r_1}^{0}   (1-e^{w})^{-\frac{1}{1+ \epsilon}}  dw = \int_{0}^{r_1}   (1-e^{-w})^{-\frac{1}{1+ \epsilon}}  dw.$$
     In  wiew of $$1-e^{-w} \geq \frac{1-e^{-r_1}}{r_1} x,  \quad  x\in  (0, r_1)$$  it results 
      $$ \int_{0}^{r_1}   (1-e^{-w})^{-\frac{1}{1+ \epsilon}}  dw\leq (1-e^{-r_1})^{-\frac{1}{1+ \epsilon}}  r_1^{\frac{1}{1+ \epsilon}}    \int_{0}^{r_1} 
        \frac{1 }{ x^{\frac{1}{1+ \epsilon}   }   } dw = \frac{1+\epsilon}{\epsilon} (1-e^{-r_1})^{-\frac{1}{1+ \epsilon}}   r_1         .$$
  %\newpage
Therefore $$[Mp\chi^{\frac{1}{p-1}}]^{\frac{p-1}{p}} \frac{T}{2} \leq  \left[   \frac{p}{(2+ \epsilon)-p(1+ \epsilon)}\left[ e^{r_0}-r_0-1 \right]^{\frac{(2+ \epsilon)-p(1+ \epsilon)}{p}  }
+ \frac{1+\epsilon}{ \epsilon}
    r_0^{ \frac{\epsilon}{1+ \epsilon}}   \right.$$
    $$ \left.+  \frac{p}{(2+ \epsilon)-p(1+ \epsilon)}  \left[ e^{-r_1}+r_1-1 \right]^{\frac{(2+ \epsilon)-p(1+ \epsilon)}{p}  }   +
\frac{1+\epsilon}{\epsilon} (1-e^{-r_1})^{-\frac{1}{1+ \epsilon}}    r_1    \right]
$$
which implies
$$\lim_{k\rightarrow 0} T=  \lim_{r_1, r_0\rightarrow 0} T =0, \qquad \mbox{ for } p\in(1,2) .$$
%$$\lim_{ r_0 \rightarrow \infty} T=0 \qquad \mbox{ for p}$$
%
To determine $\displaystyle\lim_{k\rightarrow \infty} T$ we apply the definition of $T$ 
$$T=2[Mp\chi^{\frac{1}{p-1}}]^{\frac{1-p}{p}}
\int_{-r_1}^{r_0} [ \frac{\chi}{M} k - (1+w-e^w)]^{\frac{1-p}{p}} dw.
$$
 The following bound holds
 $$\displaystyle T\geq C \int_{0}^{\frac{r_1}{2}}   [ r_1+e^{-r_1} -w-e^{-w})]^{\frac{1-p}{p}} dw
\geq C \frac{r_1}{2} \min_{w\in (0,r_1/2)}  \left[ r_1+e^{-r_1} -w-e^{-w} \right]^{\frac{1-p}{p}}    $$
 $$\qquad\qquad\qquad\qquad\qquad\qquad\qquad= C\frac{r_1}{2} \left[ \frac{r_1}{2}+e^{-r_1} - e^{-r_1/2} \right]^{\frac{1-p}{p}},    $$
and for $r_1>2$  we have that
$$  \frac{r_1}{2}+e^{-r_1} - e^{-r_1/2}  \leq \frac{r_1}{2}+1  \leq r_1.$$
Applying it in the previous inequality that verifies $T$, we get
 $$ \displaystyle T \geq c r_1^{\frac{1}{p}}$$
and then
$$\lim_{k\rightarrow \infty } T=  \lim_{r_1 \rightarrow  \infty } T = \infty, \qquad \mbox{ for } p\in(1,2) .$$
By the continuity of  $T$ with  respect to $k$,  there exist $k_n>0$ such that $T_n$,  the period of the corresponding solution $(w_n, \beta_n)$ to the energy constant $k_n$,   satisfies
$$T_n=\frac{1}{n}  $$
  and therefore
  $$w_{n}^{\prime}(0)=w_{n}^{\prime}(1)=0.$$
  We recover the solution of the original problem  $(u_n,v_n)$  defined by
  $$u_n(x)=Me^{w_n(x)}, \qquad v_n(x)= \int_{0}^{x} \beta_n (s)ds- \int_{0}^{1}\beta_n(s)ds$$
and the proof ends. \qed

\begin{remark}
We notice that,  the steady states for the 1-dimensional problem (\ref{1.1})-(\ref{3}), defined by  (\ref{3.1}), are equivalent to the solutions of the problem  
$$  \left\{  \begin{array}{ll}  \displaystyle
- u_{xx}= -  \chi  ( u w_x)_x,  &   x\in\Omega, \\[2mm]  \displaystyle
  -( | w_x|^{\frac{2-p}{p-1} }  w_x )_x = u-M, &    x\in\Omega,  \\[2mm]
 \displaystyle  w_x(0)=w_x(1)=0. &    
\end{array} \right. $$
 Since $u= Me^{\chi w}$, we have
 $$  \left\{  \begin{array}{ll}    \displaystyle
 -( | w_x|^{\frac{2-p}{p-1} }  w_x )_x= M\left(e^{\chi w}-1\right), &    x\in\Omega,  \\[2mm]
 \displaystyle w_x(0)=w_x(1)=0,&    
\end{array} \right. $$
 for $w$ satisfying  $$\int_{\Omega} e^{\chi w}=1.$$
    \end{remark}

\section*{Acknowledgments}
This work  is supported by the Project MTM2017-83391-P from MICINN (Spain).
\\

 \end{document}